\documentclass[12pt]{article}

\usepackage{amsthm,amsmath,stmaryrd,bbm,hyperref,geometry,color,amssymb}
\usepackage[all]{xy}
\usepackage{mathrsfs}

\newcommand{\po}{\left(}
\newcommand{\pf}{\right)}
\newcommand{\co}{\left[}
\newcommand{\cf}{\right]}

\newcommand{\R}{\mathbb R} 
\newcommand{\T}{\mathbb T}

\newcommand{\N}{\mathbb N} 
\newcommand{\dd}{\text{d}}

\newcommand{\new}[1]{#1}

\usepackage{graphicx}
\usepackage{enumitem,color}	

\def\pierre#1{#1}
\usepackage[normalem]{ulem}

\newtheorem{theorem}{Theorem}
\newtheorem{proposition}{Proposition}
\newtheorem{lemma}{Lemma} 

\title{Recent advances in the long-time analysis of killed degenerate processes and their particle approximation}
\author{Bertrand Cloez\footnote{MISTEA, Universit\'e Montpellier, INRAE, Institut Agro, Montpellier, France}, Lucas Journel\footnote{LJLL, Sorbonne Université, Paris},  Pierre Monmarch{\'e}\footnote{LJLL and LCT, Sorbonne Université, Paris},\\ Boris Nectoux\footnote{LMBP, Université Clermont Auvergne} and Mouad Ramil\footnote{CERMICS, Ecole des Ponts, Marne-la-Vallée, France}}

\begin{document}

%
%
%

\maketitle

\begin{abstract} We review some recent results of quantitative long-time convergence for the law of a killed Markov process conditioned to survival toward a quasi-stationary distribution, and on  the analogous question for the particle systems used in practice to sample these distributions. With respect to the existing literature, one of the novelties of these works is the degeneracy of the underlying process with respect to classical elliptic diffusion, namely it can be a non-elliptic hypoelliptic diffusion, a piecewise deterministic Markov process or an Euler numerical scheme. 

\medskip

  Nous présentons quelques résultats récents d'estimées quantitatives de convergence en temps long pour la loi de processus de Markov tués, conditionnellement à leur survie, vers une mesure quasi-stationnaire, et sur la question analogue pour les systèmes de particules utilisés en pratique pour échantillonner ces mesures. Par rapport à la littérature antérieure, l'une des nouveautés de ces travaux est le caractère dégénéré du processus de Markov sous-jacent par rapport aux diffusions elliptiques  classiques, au sens où ce peut être une  diffusion hypoelliptique non-elliptique, un processus déterministe par morceaux, ou un schéma d'Euler.  
  
  \end{abstract}
%
%

\section{Introduction}
\label{se:intro}

\subsection{Killed Markov processes}

A killed process on a space $\pierre{\mathcal D}$ is a Markov process $(X_t)_{t\geqslant 0}$ on $\pierre{\mathcal D}\cup \{\partial\}$ where $\partial \notin \pierre{\mathcal D}$ is an arbitrary cemetery point which is an absorbing point, i.e. $X_t=\partial$ for all $t\geqslant \tau =\inf\{s\geqslant 0,\ X_s=\partial\}$. In many cases of interest, it is obtained from a Markov process $(Y_t)_{t\geqslant 0}$ on a state $\pierre{\mathcal S}$ by setting $X_t=Y_t$ for $t<\tau$ and $X_t=\partial$ for $t\geqslant \tau$ where $\tau$ is defined either by a so-called soft-killing mechanism, namely, given $E$ a standard exponential random variable independent from $Y$,
\[\tau = \inf\left\{t\geqslant 0,\ E\leqslant \int_0^t \lambda(Y_s)\dd s \right\}\,,\]
with a death rate $\lambda:\pierre{\mathcal{S}}\rightarrow \R_+$, \pierre{in which case $\mathcal D=\mathcal S$,} or a so-called hard-killing mechanism, namely
\begin{equation}
\label{eq:soft}
\tau = \inf\left\{t\geqslant 0,\ Y_t\notin \mathcal D\right\}\,,
\end{equation}
where $\mathcal D$ is a given domain of $\pierre{\mathcal S}$. \pierre{Since we are only interested in the process before absorption, in these cases it is equivalent to work    with either $X$ or $Y$.  }

Some objects of interest are then, for all $t\geqslant 0$, $p_t= \mathbb P(\tau \leq t)$ the extinction probability, $\eta_t$ the law of the killed process and $\mu_t$ the law of the process conditioned to survival, given by
\begin{eqnarray}\label{def:etamu}
\eta_t(A) = \mathbb P(X_t \in A)\,,\qquad \mu_t(A) = \mathbb P(X_t \in A\ |\ \tau > t) = \frac{\eta_t(A)}{1-p_t}
\end{eqnarray}
for all measurable sets $A$ of $\pierre{\mathcal D}$. 

In most cases, extinction is almost sure, i.e. $\tau <+\infty$ almost surely. In that case, the only invariant measure of $X$ is the Dirac mass at $\partial$, which is not interesting. In fact, for killed processes, the relevant analogous to invariant measures  is given by the quasi-stationary distributions (QSD), which by definitions are the probability measures $\nu$ on $\pierre{\mathcal D}$ such that 
\[X_0 \sim \nu \quad \Rightarrow \quad \forall t\geqslant 0,\ \mu_t = \nu\,.\]
Similarly, the classical question of convergence toward equilibrium for Markov processes is replaced, for killed process, by the question of the   limit of $\mu_t$ toward a QSD as $t\rightarrow +\infty$ when the initial distribution is not a QSD.

An important property  is that, starting from a QSD, the extinction time follows an exponential distribution. In other words, if $\nu$ is a QSD, then there exists $\theta>0$ such that, if $X_0 \sim \nu$, then $p_t = 1-e^{-\theta t}$ for all $t\geqslant 0$. Equivalently, $\eta_t = e^{-\theta t} \nu$ for all $t\geqslant 0$, which means that $\nu$ is an eigenvector of the killed semi-group. In fact, when $\pierre{\mathcal D}$ is a finite set, the existence of a QSD stems from the Perron-Frobenius theorem.

\pierre{For general and extensive references on QSD, we refer the interested reader to the book \cite{CMSlivre} or the survey \cite{MV2012}.}

\subsection{Some motivations}\label{sec:motivations}

Among other things, the study of killed processes and quasi-stationary distributions is related to metastability phenomena, i.e. cases where a Markov process stays for very long times in some transient states before reaching stationarity. The convergence to a QSD within a transient state is fast and observed before the slow transitions. This arises in particular in population models and stochastic algorithms in molecular dynamics, as we describe in the following.

In population dynamics, a Markov process $(X_t)_{t\geq 0}$ is generally used to model the number of individuals. This can for example be a diffusion process on a bounded set as in the case of the evolution of gene populations in the continuous Wright-Fischer model \cite[Page 368]{ChaM16} or a continuous-time Markov chain as for a Galton-Watson process \cite[Chapter III]{AN04}. An interesting example is given by the chemostat model  \cite{CMMS13,CF22} where $X_t= (N_t, S_t) \in \mathbb{N} \times \mathbb{R}_+$ is a couple of two components: a discrete one $N_t$ and a continuous one $S_t$. 
\begin{figure}
    \centering
    \includegraphics[scale=0.5]{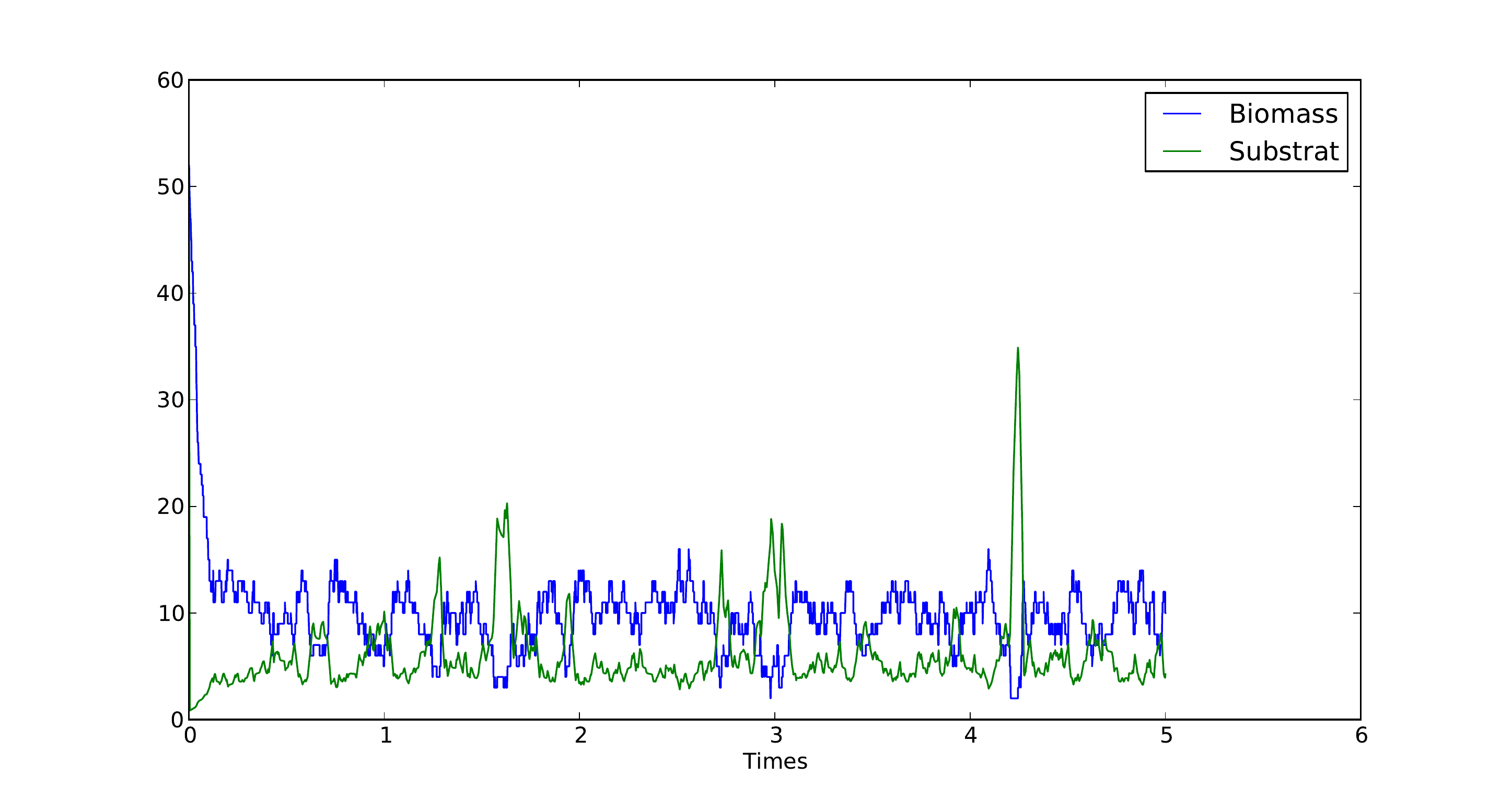}
    \caption{A sample path of the chemostat model. Here the blue curve corresponds to the discrete component $(N_t)$  and the green curve corresponds to the continuous component $(S_t)$. }
    \label{fig:chemostat}
\end{figure}
The discrete component represents the number of bacteria in a vessel and evolves as a pure jump process with the addition or removal of an individual at exponential clocks. Parameters of the jump times depend on the continuous component which evolves according an ordinary differential equation whose parameters depend on the discrete component. See  \cite{CMMS13,CF22} for details and Figure~\ref{fig:chemostat} for an example of trajectory. This process is not irreducible, belongs to a non-compact space and is hybrid, in the sense that its transition kernel does not admit a density with respect to a simple reference measure. Whatever the parameters, the population of bacteria always hits $0$ in finite time (as a consequence of a Borel-Cantelli type arguments since resources are finite), however, when the vessel is large the dynamics is close to a deterministic model having a non-trivial equilibrium. These two antagonistic behaviors, which are common for population models, should be explained by a quasi-stationary behavior.



As we shall see here, the study of QSD can also be useful in molecular dynamics methods for the sampling of the time-evolution of complex particle systems. This sampling is done using numerical schemes which are used in a wide range of application fields (biology, chemistry, materials science for instance). They allow us to compute dynamical quantities of interest for the system without resorting to the actual experiment. For instance, one would like to sample the trajectories of a system which goes from one state to another one, in order to compute the typical time to observe such transitions. This question arises for instance when we want to sample the phase transitions of a thermostated (fixed temperature) molecular system, whose evolution can be modeled by the Langevin dynamics:
\begin{equation}\label{eq:Langevin}
  \left\{
    \begin{aligned}
        &\mathrm{d}x_t=M^{-1}v_t \mathrm{d}t , \\
        &\mathrm{d}v_t=F(x_t) \mathrm{d}t -\gamma M^{-1}v_t
        \mathrm{d}t +\sqrt{2\gamma\beta^{-1}} \mathrm{d}B_t ,
    \end{aligned}
\right.  
\end{equation}
where $(x_t, v_t)\in\mathbb{R}^d\times\mathbb{R}^d$ denotes the positions and momenta of the particles at time $t\geq0$, $M$ is the mass matrix, $F$ is the interaction force (often $F$ is conservative, i.e. $F=-\nabla V$ for some potential function $V$), $\gamma > 0$ is a friction parameter, $\beta=(k_BT)^{-1}$ is proportional to the inverse temperature and $(B_t)_{t\geq0}$ is a standard Brownian motion. When the friction coefficient $\gamma$ goes to infinity, it is known that its rescaled position coordinate $(x_{\gamma t})_{t\geq0}$ converges to the solution $(\overline{x}_{t})_{t\geq0}$ of the overdamped Langevin Dynamics (see for instance~\cite[Section 2.2.4]{LelRouSto10}):
\begin{equation}\label{eq:overdamped Langevin}
    \mathrm{d}\overline{x}_t = F(\overline{x}_t)\mathrm{d} t + \sqrt{2\beta^{-1}}\mathrm{d} B_t,
\end{equation} 
  
Phase transitions are known to be often metastable, which means that the timestep used for the sampling of the dynamics is much smaller than the timescale associated with the transition event. As a practical consequence, the number of iterations needed to sample the transition event will be far too large to be done through naive simulation. The metastability comes from the existence of energetic barriers, namely the trajectory to leave the state requires to climb above a saddle point of the potential energy $V$.  

Several methods have been developed in order to circumvent this issue, among which the Parallel Replica method~\cite{Voter}. The idea behind Parallel Replica is that, as already mentioned, when the system remains trapped inside a metastable state for a sufficiently long time, its distribution becomes close to the QSD of the state. Furthermore, starting from the QSD, the first exit from the state satisfies some well-known properties which allow for the sampling in parallel of the first exit event, thus reducing the wall clock time simulation in Parallel Replica compared to standard simulation.

The main ingredients for the mathematical justification of Parallel Replica are the existence of a QSD in the state considered, as well as a long time convergence of the dynamics trapped inside the state towards the QSD. Such proofs were carried out in the literature in the case of the overdamped Langevin~\eqref{eq:overdamped Langevin} dynamics in~\cite{LebLelPer} but had yet to be extended to the Langevin dynamics~\eqref{eq:Langevin}.

\subsection{Past results and limitations }

 As far as continuous-time processes in a continuous space are concerned, there is a wide literature on quasi-stationary distributions for elliptic diffusions like the overdamped Langevin dynamics~\eqref{eq:overdamped Langevin} on smooth bounded domains. We refer for example to~\cite{GQZ,LebLelPer,V2,V3,Pinsky,CattColMelMart,CatMel}. Let us state a theorem to summarize the main   results known for the overdamped Langevin.
 
The infinitesimal generator of the overdamped Langevin dynamics is given by: 
\begin{equation}\label{L}
    \overline{\mathcal{L}}=F\cdot\nabla+\beta^{-1}\Delta, 
\end{equation}
with formal adjoint $\overline{\mathcal{L}}^*$ in $\mathrm{L}^2(\mathrm{d}x\mathrm{d}v)$ given by:  
\begin{equation*}\label{generateur adjoint overdamped}
    \overline{\mathcal{L}}^*=-\mathrm{div}(F\cdot )+\beta^{-1}\Delta. 
\end{equation*} 
The following is based on~\cite{V3,GQZ,LebLelPer,KnobPart}.
\begin{theorem}[QSD of the overdamped Langevin process]\label{qsd overdamped} Let $\beta>0$. Assume that $F\in\mathcal{C}^\infty(\mathbb{R}^{d},\mathbb{R}^{d})$ and let $\mathcal{O}$ be an open $\mathcal{C}^2$ bounded connected set of $\mathbb{R}^d$. Then there exists a unique QSD $\overline{\mu}$ on $\mathcal{O}$ for the killed process in $\mathcal{O}$, $(\overline{X}_t)_{t\geq0}$, solving~\eqref{eq:overdamped Langevin}. Furthermore, 
\begin{enumerate} 
    \item there exists $\overline{\psi}\in\mathcal{C}^2(\mathcal{O})\cap\mathcal{C}^b(\overline{\mathcal{O}})$ such that $\overline{\mu}(\mathrm{d}q)=\overline{\psi}(q)\mathrm{d}q$, where $\mathrm{d}q$ is the Lebesgue measure on $\mathbb{R}^d$,
    \item $\mathrm{Span}(\overline{\psi})$ is the eigenspace associated with the smallest eigenvalue $\overline{\theta}$ of the operator $-\overline{\mathcal{L}}^*$ with homogeneous Dirichlet boundary conditions on $\partial\mathcal{O}$, 
    \item there exist $C>0$ and $\alpha>0$ such that for all probability measures $\nu$ on $\mathcal{O}$, for all $t\geq0$, $$\big\Vert \overline{\mu}_t-\overline{\mu}\big\Vert_{TV}\leq C\mathrm{e}^{-\alpha t},$$ where $\overline{\mu}_t$ is the law of the process $(\overline{X}_t)_{t\geq0}$ conditioned to survival, as defined in~\eqref{def:etamu}, satisfying $\overline{X}_0\sim\nu$ and $\Vert\cdot\Vert_{TV}$ is the total-variation norm on the space of bounded signed measures on $\mathbb{R}^d$.
\end{enumerate}  
\end{theorem}
 
 Now, additional difficulties arise when the process is not an elliptic diffusion process, which is the case of the processes mentioned in the previous section: the kinetic Langevin dynamics~\eqref{eq:Langevin} is an hypoelliptic non-elliptic diffusion, and the chemostat model of \cite{CMMS13,CF22} is a piecewise deterministic Markov process, lacking the regularization properties of diffusions. Moreover, the study of hard-killed process with an unbounded domain (which is the case for the Langevin process since the velocity is unbounded) is also difficult. For instance even for the very simple Ornstein-Uhlenbeck process killed at the origin, there is no uniqueness of the QSD \cite{LS00}.
 
Multiple approaches were used in the literature to obtain Theorem~\ref{qsd overdamped}. To name just one, one can apply the Krein-Rutman theorem to the semigroup of~\eqref{eq:overdamped Langevin} killed outside of $\mathcal{O}$, which is compact on $\mathrm{L}^2(\mathcal{O})$, to obtain the QSD. However, such compactness property is not clear once the infinitesimal generator is not elliptic or the domain is unbounded as is the case for~\eqref{eq:Langevin} on the domain $\mathcal{D}=\mathcal{O}\times\mathbb{R}^d$. \new{Two approaches are presented in Sections~\ref{sec:Langevin 1} and~\ref{sec:Langevin 2} to establish this property in cases that cover the kinetic Langevin process \eqref{eq:Langevin}.   Section~\ref{sec:harris} presents an alternative method, based on a Doeblin-type minorization condition rather than compactness,  to get a result similar to the  Krein-Rutman theorem in non-compact space. As can be seen by comparing the general results stated in Theorems~\ref{th:cv} and \ref{thm:main}, the two methods share some similarities (a Lyapunov function is required in both cases) but each have their interest: the conditions of Theorem~\ref{thm:main} are easier to check in practice, but   the process  is required to be Feller, which is not the case for Theorem~\ref{th:cv}, while in the latter a non-trivial minorization condition has to be established. }

Finally, notice that, when it comes to sample the QSD (for instance as the first step of the Parallel Replica algorithm), a naive Monte Carlo method, consisting of running $N$ independent copies of the process up to a long time $t$ and taking the empirical distribution of the processes which have not been killed, may suffer from a high variance since the size of the sample shrinks with times (possibly, all the processes may have already died at time $t$). For this reason, systems of interacting particle are used in practice: when a process dies, another process is chosen at random uniformly and is duplicated, which maintains the size of the sample, at the cost of a small bias. The particle system is then a non-killed Markov process, which converges in long time to an invariant measure which is close, for large $N$, to the law of $N$ independent realizations of the QSD. Although many results are available in various frameworks for the convergence of the system of particles toward the killed process as $N\rightarrow \infty$, either at a fixed time $t$ or at stationary (see \cite{DM13,CC21}   and references within) obtaining quantitative convergence as $t\rightarrow \infty$ uniformly in $N$ is challenging and cannot be obtained directly from the long-time convergence of the killed process. Such a result, in the case of an elliptic diffusion, is presented in Section~\ref{sec:numeric}. The effect of the time-discretization of the  continuous-time diffusion, which is a very classical topic for the invariant measure of Markov processes but not for QSD, is also addressed, closing the gap between the theoretical results such as Theorem~\ref{qsd overdamped} and the algorithms used in practice to sample the QSD.

\section{Harris-type method for non-conservative processes}

\label{sec:harris}

\subsection{Harris theorem type assumption for Krein-Rutman theorem type conclusion}

To study $(\eta_t)$ or $(\mu_t)$ defined in \eqref{def:etamu}, it is usual to introduce the family of operators $(M_t)$ defined for all $t\geq 0$ by
$$
M_t f(x) = \mathbb{E}\left[f(X_t) \mathbf{1}_{\tau> t}  \ | \ X_0=x \right],
$$ 
where $f$ is a bounded type function and $x\in \mathcal{D}$. It is easy to show that $(M_t)$ defines a positive semigroup, but in contrast with Markov semigroup, it is not conservative in the sense that
$$
M_t \mathbf{1}_\mathcal{D} = 1-p_t \neq 1
$$
in general, where $\mathbf{1}$ is the constant function equal to $1$.

This particularity is a drawback because there is a rich base of tools to study conservative semigroups as for instance Doeblin contraction or Lyapunov methods. More precis\new{e}ly, a Markov semigroup $(P_t)$ satisfies the Doeblin minorization condition if there exists $t_0,\epsilon>0$ and a probability measure such that for all $x\in X$ 
\begin{equation}
\label{eq:doeblin}
\delta_x P_{t_0} \geq \epsilon \nu.
\end{equation}
This induces uniqueness of an invariant distribution and exponential convergence to it in total variation distance. When the state is not compact Assumption~\eqref{eq:doeblin} generally does not apply on  the whole space. Consequently, it is only supposed to hold on compact sets and an additional assumption is required to prove that the process remains in a compact set for sufficiently long times. This tightness condition is often verified through the existence of a so-called Lyapunov function: there exists a  positive function $V$, such that for a certain $t>0$,
$$
P_t V \leq C V + D,
$$
for some constants $C,D>0$ with $C<1$. This is often verified by the drift condition:
\begin{equation}
\label{eq:LyapCD}
\mathcal{L} V \leq - c V +d,
\end{equation} 
where $\mathcal{L}$ is the generator of $(P_t)$ and $c,d>0$. With this assumption, Doeblin minorization only need to hold on the sublevel set of $V$ to obtain exponential convergence to a unique invariant measure.

For conservative semigroups, these aspects are well known; see for instance \cite{MTIII,MTbook12,DMPS18,HM11}.

Several works aimed at generalizing these efficient techniques for non-conservative semi-groups. Hilbert metric and Birkhoff  contraction  yield another  powerful method for the analysis of semigroups, which has been well developed \cite{B57,Nussbaum, Seneta}. More recently, let us cite for instance \cite{KM03,KM05, MHJ21} in a general and continuous-time setting or works on discrete-time Feynman-Kac semigroups \cite{DM13,DM04}. Some works associate Doeblin-Harris techniques with Krein-Rutmann theorem and then often need strong Feller properties \cite{FRS,guillinqsd,LMR22}; See Section~\ref{sec:Langevin 1}  and \ref{sec:Langevin 2}.  Unfortunately, these results do not apply or are difficult to apply for piecewise deterministic dynamics as e.g. the chemostat model described in introduction.

In \cite{BaCG,BCGM,CG19}, some theoretical results were developed to generalize Doeblin-Harris methods for non-conservative semi-group\new{s} mainly in case of branching dynamics (i.e. in case where $M_t\mathbf{1} > \mathbf{1})$. These results are related to \cite{MM02,CMMM11,CMMS13,CV18,CV14}. A generalization of Harris Theorem is for instance described in \cite[Theorem 2.1]{BCGM} and is given below.

\begin{theorem}[Generalized Harris Theorem]
\label{th:cv}
(i) If there exist a couple of positive functions $(V,\psi)$,   $t_0>0$, $\beta >\alpha>0$, $C> 0$, $(c,d)\in(0,1]^2$, a subset $K$ and $\nu$ a probability  supported by $K$ such that $\sup_K V/\psi<\infty$ and
\begin{enumerate}[label=(A\arabic*), parsep=3mm]
\item \label{A1} $ M_{t_0}  V \leq \alpha   V + C  \mathbf{1}_K \psi$,
\item \label{A2} $M_{t_0} \psi \geq \beta \psi$,
\item \label{A3} for all positive and measurable function $f$, 
$$\inf_{x\in K} \frac{M_{t_0 }(f\psi)(x)}{M_{t_0}\psi(x)}  \geq c\, \nu(f)$$ 
\item \label{A4} for all positive integers $n$
$$ \nu\left(\frac{ M_{nt_0}\psi}{\psi}\right)\geq d \sup_{x\in K}  \frac{M_{nt_0}\psi(x)}{\psi(x)},$$ 
\end{enumerate}
then, there exists a unique triplet $(\gamma,h,\theta)$ of eigenelements of $M$ with $\gamma (h) = 1$,
{\it i.e.} satisfying for all $t\geq0$
\begin{equation}\label{vecteursales}
\gamma M_t=\mathrm{e}^{-\theta t}\gamma\qquad\text{and}\qquad M_th=\mathrm{e}^{-\theta t}h.
\end{equation}
Moreover,  there exist $C,\omega>0$ such that for all $t\geq0$ and measure $\mu$, 
\begin{equation}\label{eq:conv_norm}
 \sup_{\left|f/V \leq 1\right| \leq 1} \left| e^{\theta t}\mu M_t f-\mu(h)\gamma(f)\right| \leq C \mu(V) e^{-\omega t}.
\end{equation}
(ii) Assume that there exist a positive measurable function $V,$ a triplet $(\gamma,h,\theta)$ and constants $C,\omega>0$ such that
\eqref{vecteursales} and \eqref{eq:conv_norm} hold. Then, the couple $(V,h)$ satisfies Assumptions \ref{A1},\ref{A2},\ref{A3}, and \ref{A4}.
\end{theorem}

By differentiating \eqref{vecteursales}, the triplet $(\gamma,h,\theta)$ is a triplet of eigenelements for the infinitesimal generator $\mathcal{A}$ of $(M_t)_{t\geq0},$
that is $\gamma\mathcal A=-\theta\gamma$ and $\mathcal A h=-\theta h$ where the (unbounded) operator $\mathcal A$ is defined by $\mathcal A=\lim_{t\to0}\frac1t(M_t-I).$

The assumptions of Theorem~\ref{th:cv} hold on compact space when $(M_t)$ admits a (upper and lower) bounded density with respect to a fixed measure $\nu$. It is known since the 60's, by Birkhoff \cite{B57} for instance, that this type of regularity hypothesis implies the existence of eigen-elements and the convergence. Unfortunately, assuming bounded density does not include simple examples such as the chemostat model.

Assumption~\ref{A4} overcomes this problem and \new{allows to treat such an} example. This type of condition seems to appear for the first time in \cite{MM02} without the Lyapunov condition. Article \cite{MM02} only states a contraction inequality, but this is the main step of the proof for such results. In \cite{CV14}, Champagnat and Villemonais highlight these conditions and prove that they are in fact equivalent to the exponential convergence \eqref{eq:conv_norm}. They deduce the existence of eigenelements in the context of absorbed Markov processes (sub-conservative semigroups). In \cite{BaCG}, this theorem is extended for general inhomogeneous-time and non-conservative semi-groups. It includes for instance periodic and asymptotically homogeneous ones.

As explained before, requiring a Doeblin type minorization assumption on the whole space is often too strong when the state space is not compact, Assumptions \ref{A1} and \ref{A2} then relax this constraint. This Lyapunov assumption seems to appear for the first time in \cite[Theorem 4.2]{CMMM11} and in the proof of \cite[Theorem 4.1]{CMMS13}. However, combining this assumption with the assumptions \ref{A3} and \ref{A4} to obtain uniqueness of eigenelements and exponential convergence comes from \cite{CV18}. This significant result has been demonstrated in the framework of sub-conservative semigroups. Result \cite[Theorem 2.1]{BCGM} generalizes this result for non-conservative semi-groups, and non-bounded functions $\psi$, and prove that these assumptions are in fact necessary. Let us also point out the follow-up paper \cite{CV20} of \cite{BCGM,CV18} on these aspects. However, Assumption \ref{A4} of Theorem~\ref{th:cv} seems simpler to verify that the analogous of \cite{CV20} and the proof of \cite{BCGM} leads some quantitative estimates. More precisely, we have that $\gamma(V)<+ \infty$, $\gamma \gg \nu$ (in the Radon-Nikodym sense),
$$
\beta \leq e^{-\theta t_0} \leq \alpha + C,  
$$
and 
$$
c_1 \left( \frac{\psi}{V} \right)^q \psi \leq h \leq c_2 V
$$
where $c_1,c_2>0$, $q\in (0,1)$ and the others constants come from Theorem~\ref{th:cv}.

To conclude this section let us mention that from \cite{CG19}, Condition~\ref{A4} is verified when there exist $t_0, \epsilon>0$ such that
\begin{equation}
\label{eq:A4irr}
\inf_{x,y \in K} \mathbb{P}\left(\exists s \leq t_0, X_s=y \ | \ X_0=x \right) \geq \epsilon.
\end{equation}

This condition will be verified in almost all examples of Subsection~\ref{se:exe}

\subsection{Some (very) simple examples}
\label{se:exe}
In this section, we detail several examples where one of the assumption\new{s} does not hold in order to highlight their role.

Note however that all assumptions are naturally connected and that one can postpone the difficulty of an  assumption by changing the compact $K$ or the Lyapunov functions for instance.

\subsubsection{Aperiodicity and Assumption~\ref{A3}}

\begin{figure}
    \centering
    \includegraphics[scale=0.4]{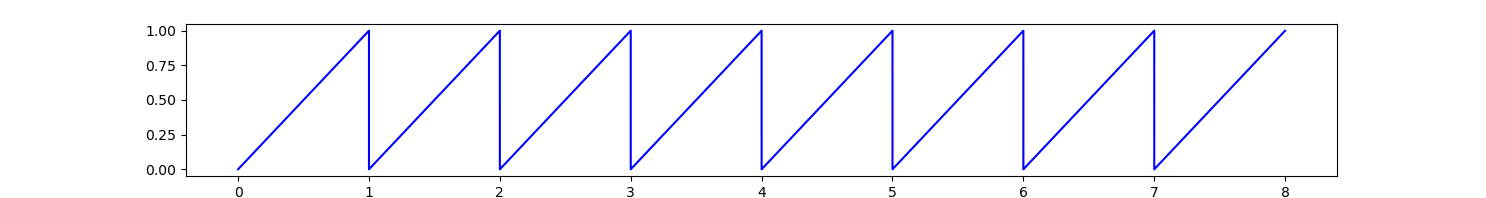}
    \caption{A sample path of the process described in Equation~\ref{eq:per}. }
    \label{fig:per}
\end{figure}

In addition to the irreducibly assumption, Assumption \ref{A3} hides an aperiodicity assumption. Indeed the (unabsorbed) Markov process $(X_t)_{t\geq 0}$ on $[0,1)$ defined, for all $t\geq 0$, by
\begin{equation}
\label{eq:per}
X_t= X_0 + t \text{ mod } 1 
\end{equation}
verifies all the assumptions except \ref{A3}. Even if the trajectories pass into all the points of $[0,1)$, it does not pass at the same time $t_0$ uniformly over its stating point. See Figure~\ref{fig:per}.

\begin{figure}
    \centering
    \includegraphics[width=14cm, height=5cm]{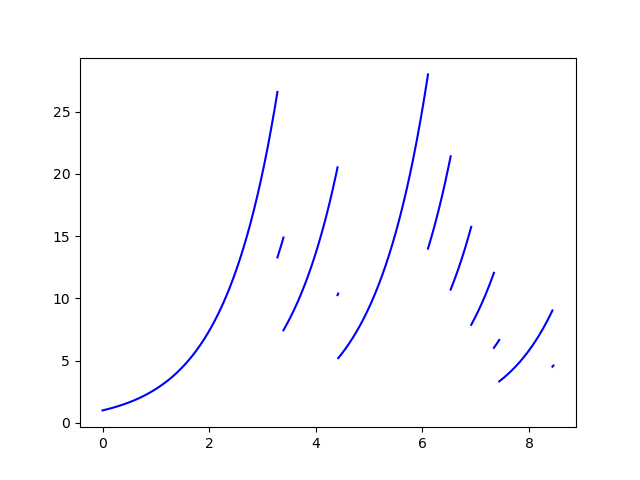}
    \caption{A sample path of the process generated by \eqref{eq:gf}. }
    \label{fig:gf}
\end{figure}

A less simple and more interesting example is given by the growth-fragmentation process studied in \cite{BDG,GM19} (in a different formalism). This process grows exponentially at a certain rate $\alpha$ and has multiplicative jumps from $x$ to $r x$, where $r\in (0,1)$ is a fixed deterministic value, at a rate $B(x)$ depending on its state. We can kill the process at some rate $\lambda$ as defined in Equation~\ref{eq:soft}. Its generator reads
\begin{equation}
\label{eq:gf}
    \mathcal{A} f(x)= \alpha x f'(x) + B(x)(f(r x) -f(x)) + \lambda(x) f(x),
\end{equation}
for smooth functions $f$ and $x>0$. See Figure~\ref{fig:gf}

In contrast with the previous trivial deterministic process, this growth-fragmentation has a random dynamics but its law is concentrated in fixed time in a discrete space (when the starting distribution is a Dirac mass). It then can not converge to its unique QSD, which has a Lebesgue density.

\subsubsection{Probability of survival and Assumption~\ref{A4}}
\label{se:2pt}
Let us consider the Markov process which is maybe the simplest one. That is, we consider that $\mathcal{D}= \{1,2\}$ and $X_t$ jumps from $2$ to $1$ at rate $a$ and cannot jump from $1$ to $2$. Similarly, it dies from $1$ at rate $b$ and cannot die from $1$. It is a pure death process whose transition are represented in Figure~\ref{fig:markov}.
Its study is trivial and based on the spectral property of a $2\times 2$ matrix.
We have three cases :

\begin{figure}
    \centering
    \includegraphics{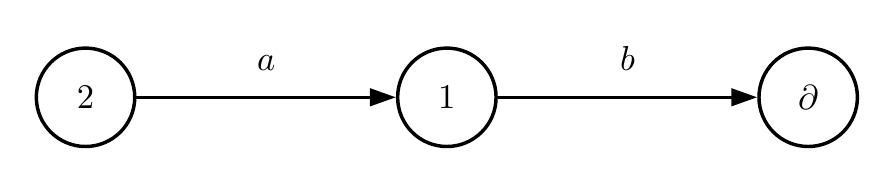}
    \caption{The transition diagram of the Markov chain described in Subsection~\ref{se:2pt} }
    \label{fig:markov}
\end{figure}

\begin{itemize}
\item $b<a$ : We have exponential convergence to the unique QSD $\delta_2$.
\item $b=a$ : We have a slow convergence to the unique QSD $\delta_2$.
\item $b>a$ : We have two QSD:
$$
\frac{a}{b} \delta_1 + \frac{b-a}{b}\delta_2, \qquad \delta_2,
$$
\new{and, furthermore,} the convergence depends on the initial condition : if $\mu_0(\{2\})=\eta(\{2\})>0$ then the conditional law $\mu_t$ converges to $\frac{a}{b} \delta_1 + \frac{b-a}{b}\delta_2$.
\end{itemize}

We see in the third case that $ \delta_2 M_t\mathbf{1}\sim e^{-at}$ and $ \delta_1 M_t\mathbf{1} \sim e^{-b t}$  are not of the same   order and consequently \ref{A4} fails and then the conclusion of Theorem~\ref{th:cv} doesn't hold (which is what we naturally and simply observe here).

This process is trivial to study but it is an interesting example to understand the algorithms which approximate the quasi-stationary distribution. For instance the limit in large number of particles and in time of the well used Fleming-Viot algorithm do not commute for this example. Some question on its approximation remains open see \cite{BCP}. Commutation of limits for QSD approximation, as detailed in Section~\ref{sec:numeric} below is then of prime importance.

A less simple and more interesting example is given by the house-of-card model studied in \cite{CG22} and references therein. This process belongs to $[0,1]$ and is redrawn uniformly at random over this interval at rate $1$. It is soft-killed at rate $x\mapsto c x^q$, for some $c,q\geq 0$. We can show the following :
 
 \begin{itemize}
\item If $q> 1 - 1/c$ then we have exponential convergence to a unique explicit QSD (which has a bounded Lebesgue density). This in particular always verified when $q\geq 1$ or $c<1$. The convergence can be established in various distances with explicit rates.
\item If $c=q+1$ then we have two QSD : one, denoted by $\nu$, has a Lebesgue density proportional to $x\mapsto x^{-q}$ and the second is $\delta_0$. When $\mu_0(\{0\})=0$, then we have a (explicit) slow convergence of the conditional laws $\mu_t$ to the QSD $\nu$ (which has a Lebesgue density proportional to $x\mapsto x^{-q}$). When $\mu_0(\{0\})>0$ then $\mu_t \to \delta_0$. Moreover, if $q <1/2$ then $\eta_t e^{\theta t} \to c_0 \nu$, for some $c_0>0$ although if $q\geq 1/2$ then $\eta_t e^{\theta t} \to 0$.
\item If $c< q+1$ then the conditional law converges, in mean,  to the unique degenerate QSD:
$$
\nu(dx) = \left(1 - \frac{c}{q+1}\right) \delta_0(dx) + \frac{x^{-q}}{c} dx
$$
\end{itemize}

For details, additional results and references, see \cite{CG22}. Note that, in this example, Assumption~\eqref{eq:A4irr} fails.

\subsubsection{Lyapunov function and Assumption~\ref{A1} and \ref{A2}}

The Lyapunov assumptions require to find a kind of sub-solution and super-solution to the eigenvalue problem where the associated sub-eigenvalue $\beta$ is strictly larger than the super-eigenvalue $\alpha$. Let us show here a well-studied example of Markov processes where we easily find two Lyapunov functions in such a way that the condition $\alpha=\beta$ become\new{s} a critical situation.

Let us consider the birth and death process on $\mathbb{N}^*$ generated by
$$
\mathcal{A}f(n) = b(f(n+1)-f(n)) + d(f(n-1)-f(n)),
$$
for every $n\geq 2$ and bounded function $f$, with $b>d$, and
$$
\mathcal{A}f(1) = b_1(f(2)-f(1)) - d_1 f(1),
$$
with some $b_1,d_1>0$.  It is known from \cite{vD03} that $e^{\theta t} \eta_t$ converge to a unique QSD if and only if
$$
(\sqrt{b} -\sqrt{d} )^2 + b_1 \left( \sqrt{d/b} - 1 \right) - d_1 >0.
$$
This result can easily be recovered (and extended) by using Theorem~\ref{th:cv} with $V$ and $\psi$ of the form $q\mapsto q^n$ with appropriate $q>0$. See \cite{BCGM} for details. Note that \new{the condition~\eqref{eq:A4irr}} is trivially satisfied.

\section{The Langevin process (Round 1)}\label{sec:Langevin 1}
Let $F\in\mathcal{C}^\infty(\mathbb{R}^{d},\mathbb{R}^{d})$, $\gamma\in\mathbb{R}$, $\sigma>0$. In this section, we shall consider the following Langevin dynamics:
\begin{equation}\label{eq:Langevin 1}
  \left\{
    \begin{aligned}
        &\mathrm{d}x_t=v_t \mathrm{d}t , \\
        &\mathrm{d}v_t=F(x_t) \mathrm{d}t -\gamma v_t
        \mathrm{d}t +\sigma \mathrm{d}B_t ,
    \end{aligned}
\right.  
\end{equation}
Its infinitesimal generator is given for $(x,v)\in\mathbb{R}^d\times\mathbb{R}^d$ by:
\begin{equation} 
    \mathcal{L} =  v\cdot\nabla_x+F(x)\cdot\nabla_v -\gamma  v\cdot\nabla_v+ \frac{\sigma^2}{2}\Delta_v, 
\end{equation}
with formal adjoint $\mathcal{L}^*$ in $\mathrm{L}^2(\mathrm{d}x\mathrm{d}v)$:
\begin{equation} 
    \mathcal{L}^*=-v\cdot\nabla_x -F(x)\cdot\nabla_v+\gamma \mathrm{div}_v(v  \cdot )+\frac{\sigma^2}{2}\Delta_v. 
\end{equation}

The metastable states often correspond to the basins of attraction of a potential $V$ when the force $F$ in~\eqref{eq:Langevin} is conservative, i.e. $F=-\nabla V$. Therefore, we shall consider metastable states written as cylinders in the position and momenta coordinates as follows:
\begin{equation}\label{eq:def domain D}
    \mathcal D=\mathcal{O}\times\mathbb{R}^d,
\end{equation}
where $\mathcal{O}$ is an open $\mathcal{C}^2$ bounded connected set of $\mathbb{R}^d$.   
\subsection{Quasi-stationary distribution: existence and long-time convergence}
The justification of the Parallel Replica method mentioned in Section~\ref{sec:motivations} for the Langevin dynamics~\eqref{eq:Langevin 1} in $\mathcal{D}$ first requires the existence of a QSD on $\mathcal{D}$ for the associated process $(X_t)_{t\geq0}$ killed outside of $\mathcal{D}$. Such existence can be obtained through the application of the Krein-Rutman theorem to the semigroup of $(X_t)_{t\geq0}$, which shall provide a left eigenvector of this semigroup corresponding to the QSD density. To apply the Krein-Rutman theorem, this semigroup is required to be compact in a Banach space. To prove the compactness the idea of the proof is to first show that the semigroup is an integral Hilbert-Schmidt operator in $\mathrm{L}^2(\mathcal D)$ hence compact in $\mathrm{L}^2(\mathcal D)$ using appropriate Gaussian upper-bounds for its transition density obtained in~\cite[Theorem 2.19]{LelRamRey}. Second, this compactness is propagated across the Banach space $\mathcal{C}^{b}(\overline{\mathcal D})$ of continuous and bounded functions on $\overline{\mathcal D}$ using the continuity of the semigroup, leading to the following result in~\cite[Theorem 2.9]{LelRamRey2}. 

\begin{theorem} [Compactness of the killed semigroup]\label{compactness} For any $t \geq0$, $p \in [1,+\infty]$ and $f \in \mathrm{L}^p(\mathcal D)$, the quantity
\begin{equation}\label{def semigroupe}
  M_t f : (x,v) \in \overline{\mathcal D} \mapsto \mathbb{E}\left[f(X_t)\mathsf{1}_{\tau>t}\vert X_0=(x,v)\right]
\end{equation} 
is well-defined. Besides, 
\begin{enumerate} 
  \item the family of operators $(M_t)_{t \geq 0}$ is a semigroup on $\mathrm{L}^p(\mathcal D)$ and on $\mathcal{C}^{b}(\overline{\mathcal D})$. 
  \item For any $t>0$, the operator $M_t$ is compact from $\mathrm{L}^p(\mathcal D)$ to $\mathrm{L}^p(\mathcal D)$, and from $\mathcal{C}^{b}(\overline{\mathcal D})$ to $\mathcal{C}^{b}(\overline{\mathcal D})$. 
\end{enumerate} 
\end{theorem}

The application of the Krein-Rutman theorem to the operator $M_t$ for $t>0$  on the Banach space $\mathcal{C}^{b}(\overline{\mathcal D})$ then yields the existence of a QSD on $\mathcal{\mathcal D}$ of $(X_t)_{t\geq0}$, with a density belonging to $\mathcal{C}^{b}(\overline{\mathcal D})$. The uniqueness of this QSD follows from the positivity of the transition density of the semigroup using a Harnack inequality provided in~\cite{LelRamRey}. In addition, one can also obtain a spectral interpretation of this QSD using Itô's formula, see the proofs of~\cite[Theorems 2.14, 2.17]{LelRamRey2} for more details.

\begin{theorem}[Existence and uniqueness of a QSD]\label{thm:spectral interpretation} There exists a unique QSD $\mu$ on $\mathcal D$ for the killed process $(X_t)_{t\geq0}$. In addition, there exists a smooth positive function $\psi$ such that
$$\mathrm{d}\mu(x,v)=\psi(x,v)\mathrm{d}x\mathrm{d}v.$$
Moreover, $\psi$ is the unique solution in $\mathcal{C}^2(\mathcal D)\cap\mathcal{C}^b(\overline{\mathcal D})$ to the eigenvalue problem
\begin{equation*} 
  \left\{
    \begin{aligned}
        \mathcal{L}^*\psi(x,v)  &=-\theta \psi(x,v)&&\quad (x,v)\in \mathcal D,\\
        \psi(x,v) &=0&&\quad x\in\partial\mathcal{O}, v\cdot n(x)<0,
    \end{aligned}
\right.  
\end{equation*} 
where $n(x)$ is the unitary outward normal vector at $x\in\partial\mathcal{O}$. 
\end{theorem} 

In addition, the compactness of the killed semigroup allows to obtain a spectral decomposition of the killed  semigroup. This gives in particular the long-time asymptotics of the law of the dynamics conditioned to remain in $\mathcal D$, see~\cite[Theorem 2.22]{LelRamRey2}. 
\begin{theorem}[Convergence to the QSD in total variation] There exists $\alpha>0$ such that for any probability measure $\nu$ on $\mathcal D$, there exists $C_\nu>0$ such that for all $t\geq0$, 
\begin{equation}\label{eq:cv of conditional distrib}
    \left\|\mu_t- \mu\right\|_{TV} \leq C_\nu \mathrm{e}^{-\alpha t},
\end{equation}
where $\mu_t$ is the law of $(X_t)_{t\geq0}$ conditioned to survival defined in~\eqref{def:etamu} satisfying $X_0\sim\nu$, and $\|\cdot\|_{TV}$ denotes the total-variation norm on the space of bounded signed measures on $\mathbb{R}^{2d}$.
\end{theorem} 
The QSD on $\mathcal D$ can thus be seen as the longtime limit of the dynamics conditioned to stay inside $\mathcal D$. This proposition is useful to understand what is a metastable state. A metastable state is a state such that the exit event from this state happens on a timescale larger than the local equilibration time, namely the time to observe the convergence to the QSD in~\eqref{eq:cv of conditional distrib}.
  
Furthermore, if the Langevin dynamics~\eqref{eq:Langevin 1} is initially distributed according to the QSD in $\mathcal D$, then we can explicitly provide the law of the first exit event from $\mathcal D$, which is given by the couple of the first exit time from $D$ and its associated exit point. As mentioned in the introduction, recall that the first exit time from a domain when starting from a QSD is well known to be distributed according to an exponential law; this is a general property for Markov processes. By contrast, the law of the first exit point is strongly dependent on the dynamics. In the case of the Langevin dynamics~\eqref{eq:Langevin 1}, this is the first result to compute the law of the first exit point starting from the QSD.

\begin{proposition}\label{prop:first exit event law}
Let us assume that $X_0$ is distributed according to the QSD $\mu$ in $\mathcal D$. Then the law of the couple $(\tau,X_{\tau})$ is fully characterized by the following properties:
\begin{itemize}
    \item $\tau$ follows the exponential law of parameter $\theta$ where $\theta$ is defined in Theorem~\ref{thm:spectral interpretation};
    \item $\tau$ is independent of $X_{\tau}$;
    \item The law of $X_{\tau}$ is given by: for any bounded function $f:\partial \mathcal D\mapsto\mathbb{R}$, 
    $$\mathbb{E}\left[f(X_{\tau})\vert X_0\sim\mu \right]=\frac{1}{\theta}\int_{\partial \mathcal D}\left\vert v\cdot n(x)\right\vert\psi(x,v)f(x,v)\sigma_{\partial \mathcal{O}}(\mathrm{d}x)\mathrm{d}v,$$
    where $n(x)$ is the unitary outward normal vector at $x\in\partial\mathcal{O}$ and $\sigma_{\partial \mathcal{O}}$ denotes the Lebesgue measure on the surface $\partial\mathcal{O}$.
\end{itemize}
\end{proposition} 
As a result, the existence of a QSD and the long-time convergence to the QSD of the Langevin dynamics trapped in $\mathcal D$ ensure that the Parallel Replica method also applies for the Langevin process trapped in $\mathcal{D}$ following the reasoning made in~\cite[Proposition 5]{LebLelPer}. 
 
\subsection{Quasi-stationary distribution: overdamped limit}

In this section we shall assume that $\gamma>0$ and $\sigma=\sqrt{2\gamma\beta^{-1}}$ for some $\beta>0$ in~\eqref{eq:Langevin 1}. The aim of this section is to describe the behavior of the Langevin dynamics~\eqref{eq:Langevin 1} when the friction coefficient $\gamma$ goes to infinity. This will allow us in particular to explicit the overdamped limit of the QSD. 

It is well known in the literature that the law of the rescaled position of the Langevin dynamics converges to the law of the overdamped Langevin dynamics~\eqref{eq:overdamped Langevin}, i.e. for all $T>0$:
\begin{equation}\label{eq:OD limit}
    \mathrm{Law}((x_{\gamma t})_{t\in[0,T]})\underset{\gamma\rightarrow\infty}{\longrightarrow}\mathrm{Law}((\overline{x}_{t})_{t\in[0,T]}).
\end{equation}
However, this result does not provide an overdamped limit for the couple of the position and velocity vectors. Such a limit was provided recently in~\cite{Ram} where the following limit is obtained: 
\begin{equation}\label{eq:OD limit couple}
    \mathrm{Law}((x_{\gamma t})_{t\in[0,T]},v_{\gamma T})\underset{\gamma\rightarrow\infty}{\longrightarrow}\mathrm{Law}((\overline{x}_{t})_{t\in[0,T]},Z),
\end{equation}
where $Z  \sim  \mathcal{N}_{d}(0,\beta^{-1} I_d)$ is a Gaussian vector independent of the process $(\overline{x}_{t})_{t\in[0,T]}$. The proof of~\eqref{eq:OD limit couple} consists in perturbing the Brownian noise of the Langevin dynamics~\eqref{eq:Langevin 1} by a vanishing term when $\gamma$ goes to infinity such that the position and velocity become independent. Then it remains to consider the overdamped limit of the position and velocity coordinates separately, which is much easier. In addition, since the perturbation converges to zero when $\gamma$ goes to infinity, the overdamped limit of the perturbed process is the same as the overdamped limit of the Langevin dynamics.   
Following this result one can identify explicitly the overdamped limit of the QSD for the Langevin dynamics.

\begin{proposition}
Let $\mu^{(\gamma)}$ (resp. $\overline{\mu}$) be the QSD of $(X_t)_{t\geq0}$ on $D=\mathcal{O}\times\mathbb{R}^d$ (resp. $(\overline{X}_t)_{t\geq0}$ on $\mathcal{O}$)  satisfying~\eqref{eq:Langevin 1} (resp.~\eqref{eq:overdamped Langevin}).  Then, 
\begin{equation}\label{eq:cv overdamped qsd}
    \mu^{(\gamma)}(\mathrm{d}x\mathrm{d}v)\underset{\gamma\rightarrow\infty}{\longrightarrow}\overline{\mu}(\mathrm{d}x)\frac{\mathrm{e}^{-\frac{\beta\vert v\vert^2}{2}}}{(2\pi\beta^{-1})^{\frac{d}{2}}} \mathrm{d}v,
\end{equation}
in terms of weak convergence measures.  
\end{proposition}  
In particular, if one considers the marginal in position of the Langevin QSD then it converges weakly to the QSD of the overdamped dynamics on $\mathcal{O}$ as a consequence of the above proposition.


 \section{Langevin processes (round 2)}
\label{sec:Langevin 2}

In this section, we present  the  main results of~\cite{guillinqsd,guillinqsd2}. These results first aim at giving necessary conditions   for  existence and  uniqueness of a quasi-stationary distribution $\mu_{\mathcal D}$ of
a strongly Feller Markov process $(X_t,t\ge 0)$ killed when it exits a subdomain $\mathcal  D$, see Section \ref{sec.Secgeneral}. 
These results are then applied to Langevin processes with continuous or singular potentials, see Sections \ref{sec.Sec-appl1} and \ref{sec.Sec-appl2}.
 

\subsection{A general result}
\label{sec.Secgeneral}

\subsubsection{Notations and assumptions}
Let $(X_t,t\ge 0)$ be a time homogeneous Markov process valued in a  Polish  space $\mathcal S$, with c\`adl\`ag paths and satisfying the strong Markov property, which is defined on the filtered probability space $(\Omega, \mathcal F, (\mathcal F_t)_{t\ge 0}, (\mathbb P_{\mathsf x})_{x\in \mathcal S})$ where $\mathbb P_{\mathsf x}(X_0=\mathsf x)=1$ for all $x\in \mathcal S$. Its transition probability semigroup is denoted by $(P_t,t\ge 0)$. 

Let $\mathcal B(\mathcal S)$ be the Borel $\sigma$-algebra of $\mathcal S$, $b\mathcal B(\mathcal S)$ the space of all bounded and Borel measurable functions $f$ on $\mathcal S$ (its norm will be denoted by $f\in b\mathcal B(\mathcal S)\mapsto \Vert f \Vert_{b \mathcal B(\mathcal S)}=\sup_{\mathsf x\in \mathcal S}\vert f(\mathsf x)\vert$).  The space $\mathcal D([0,T],\mathcal S)$ of $\mathcal S$-valued c\`adl\`ag paths   defined on $[0,T]$ is equipped with the Skorokhod topology.

We suppose that

\begin{enumerate}
\item[\textbf{(C1)}] \textbf{(strong Feller property)} There exists $t_0>0$ such that for each $t\ge t_0$, $P_t$ is strong Feller, i.e. $P_tf$ is continuous on $\mathcal S$ for any $f\in b\mathcal B(\mathcal S)$.

\item[ \textbf{(C2)}]  \textbf{(trajectory Feller property)} For every $T>0$, $\mathsf x\to\mathbb P_{\mathsf x}(X_{[0,T]}\in \cdot)$  (the law of $X_{[0,T]}:=(X_t)_{t\in [0,T]}$) is continuous from $\mathcal S$ to the space $\mathcal M_1(\mathcal D([0,T], \mathcal S))$ of probability measures on  $\mathcal D([0,T],\mathcal S)$,
equipped with the weak convergence topology.
\item[ \textbf{(C3)}] \textbf{(Lyapunov function condition)} There exist a continuous  function $\mathsf W : \mathcal S\to [1,+\infty[$, belonging   to the extended domain $\mathcal D_e(\mathcal L)$ of the generator of $(X_t,t\ge 0)$,   two sequences of positive constants $(r_n)$ and $(b_n)$ where $r_n\to +\infty$, and an increasing sequence of compact subsets $(K_n)$ of $\mathcal S$, such that
$$
-\mathcal L \mathsf W ({\mathsf x}) \ge r_n \mathsf W({\mathsf x}) - b_n \mathsf 1_{K_n}({\mathsf x}), \ q.e,
$$
where $\mathsf 1_{K_n}$ is the indicator function of $K_n$. 
\end{enumerate}

Now let $\mathcal  D$ be a non-empty open domain of $\mathcal S$, different from $\mathcal S$. Consider the first exit time of $\mathcal  D$
$$
\tau:=\inf\{t\ge0, X_t\in \mathcal D^c\}
$$
where $\mathcal D^c= \mathcal S\backslash \mathcal D$ is the complement of $\mathcal  D$.
The transition semigroup of the killed process $(X_t, 0\le  t< \tau)$ is: 
\begin{equation}\label{killedsg}
M_t f(\mathsf x) = \mathbb E_{\mathsf x} [\mathsf 1_{t<\tau} f(X_t)], \ \  \text{$t\ge 0$, $\mathsf x\in \mathcal D$, $f\in b\mathcal B(\mathcal D)$}.
\end{equation} 
Assume the following on the killed process $(X_t, 0\le  t< \tau)$:
 
\begin{enumerate}
\item[\textbf{(C4)}] For $t\ge 0$, $M_t$ is  Feller, i.e. if $f$ is bounded and continuous over $\mathcal D$, so is $M_tf$. 

\item[\textbf{(C5)}] There exists $t_0>0$ such that for all $t\ge t_0$,  ${\mathsf x}\in \mathcal D$, and non-empty open subsets $O$ of $\mathcal  D$,
$$
M_t({\mathsf x},O)>0
$$
(we can assume this $t_0>0$ is the same as the one in \textbf{(C1)}), 
and there is some $\mathsf x_0\in \mathcal D$ such that $\mathbb P_{\mathsf x_0}(\tau<+\infty)>0$.
\end{enumerate}

\subsubsection{The general result for existence and uniqueness of the quasi-stationary distribution}

In this section, we provide  a slightly less general result than~\cite[Theorem 2.2]{guillinqsd}.  Indeed, Theorem~\ref{thm:main} below is actually still valid when replacing \textbf{(C4)} by the less stringent assumption that for $t\ge 0$, $M_t$ is weakly Feller.

\begin{theorem}
\label{thm:main}
Assume that \textbf{(C1)}$\to$\textbf{(C5)} hold. \new{Take any }  $p\in]1,+\infty[$. Then, there exists only one quasi-stationary distribution  $\mu_{\mathcal D}$  of the process $(X_t,t\ge 0)$ in $\mathcal D$  satisfying
$$\mu_{\mathcal D}(\mathsf W^{1/p}):=\int_{\mathcal D} \mathsf W^{1/p}({\mathsf x}) \mu_{\mathcal D}(d{\mathsf x})<+\infty.$$
In addition, there are some constants $\delta>0$ and $C\ge 1$ such that for any initial distribution $\nu$ on $\mathcal  D$ with $\nu(\mathsf W^{1/p})<+\infty$,
\begin{equation}\label{eq.thm-maina}
\big |\mathbb P_\nu(X_t\in A| t<\tau)-\mu_{\mathcal D}(A)\big |\le C e^{-\delta t} \frac{\nu(\mathsf W^{1/p})}{\nu(\varphi)}, \ \forall A\in\mathcal B(\mathcal D), t>0,
\end{equation}
where $\varphi:\mathcal D\to \mathbb R$ is continuous, $\varphi/ \mathsf W^{1/p}$ is bounded, and $\varphi$ is positive  within  $\mathcal  D$. 
\end{theorem}

Let us give some remarks on this theorem. 
If $\mathsf W$ is bounded  over $\mathcal D$, the quasi-stationary distribution inside $\mathcal D$ is unique. When $\mathcal S$ is moreover compact, Theorem~\ref{thm:main} is   valid without assuming $\textbf{(C3)}$ and with $\mathsf W=1$ there (i.e. there exists a unique quasi-stationary distribution in $\mathcal M_1(\mathcal D)$ and \eqref{eq.thm-maina} holds for any $\nu \in \mathcal M_1(\mathcal D)$). \new{On the contrary, when $\mathsf W$ is unbounded, Theorem~\ref{thm:main} states that there is a unique QSD which integrate $\mathsf W^{1/p}$ for all $p>1$ (in the statement the QSD $\mu_{\mathcal D}$ a priori depends on $p$, but then, clearly, it doesn't, since $\mu_{\mathcal D}(\mathsf W^{1/p}) < \infty$ for some $p>1$ implies $\mu_{\mathcal D}(\mathsf W^{1/q}) < \infty$ for all $q\geqslant p$), but it doesn't prevent the existence of other QSD $\nu$ with $\nu(\mathsf W^{1/p})=+\infty$ for all $p>1$. This is for instance the case for the Ornstein-Uhlenbeck process killed at $0$ (see \cite{LS00}), to which Theorem~\ref{thm:main} applies with $\mathsf W(\mathsf x) = e^{\alpha \mathsf x^2}$ with $\alpha<1/2$.  }

The intuition behind the assumptions  \textbf{(C1)}$\to$\textbf{(C5)}  is the following.   
It is now well known 
 that a Lyapunov  condition provides a spectral gap for a non killed semigroup in some weighted Banach spaces~\cite{DMT1995,MT1993}. Using results of \cite{Wu2004} characterizing the essential spectral radius of a positive operator (combined with the regularity conditions \textbf{(C1)} and \textbf{(C2)} on the non killed process $(X_t,t\ge 0)$),   the Lyapunov type condition  \textbf{(C3)}   is used    to  get  a spectral gap for $M_t$ in    $b_{\mathsf W^{1/p}}\mathcal B(\mathcal D)$ (the space of  measurable functions $f$ over $\mathcal D$ such that $f/\mathsf W^{1/p}$ is bounded). It is also known that  an irreducibility condition and a regularity condition on a non killed semigroup ensure  the uniqueness of the invariant measure. It is similar for killed processes:   we show that   the regularity condition \textbf{(C4)}   on $M_t$ and  the irreducibility condition \textbf{(C5)} imply that $M_t$ admits  a principal eigenvalue.

Let us also refer to~\cite{CV14,CV18,CattColMelMart,kolb2012quasilimiting,BCGM,ferre2,velleret2018unique} and references therein for other general criteria for existence and uniqueness of a quasi-stationary distributions.


\subsection{Langevin processes with continuous coefficients}
\label{sec.Sec-appl1}
 

In this section we apply Theorem \ref{thm:main} to Langevin processes   with continuous coefficients (also called hypoelliptic damped Hamiltonian systems). Let us introduce these processes. 
Let  $d\ge 1$ and     $(\Omega, \mathcal F, (\mathcal F_t)_{t\ge 0}, \mathbb P)$ be a  filtered probability space.  Let   $( X _t=(x_t,v_t), t\ge 0  )$ be the solution of the following hypoelliptic stochastic differential equation on $\mathbb R^{2d}$:
 \begin{equation} \label{eq.hypo}
\left\{
\begin{aligned}
 d  x_t &= v_t dt,\\
dv_t&=-\nabla V(x_t)dt - c(x_t,v_t) v_t dt +  \Sigma(x_t,v_t)  dB_t,
\end{aligned}
\right.
\end{equation}
  where $ ( B _t, t\ge 0  )$ is a standard $d$-dimensional Brownian motion on    $(\Omega, \mathcal F, (\mathcal F_t)_{t\ge 0}, \mathbb P)$. Here the state space is $\mathcal S=\mathbb R^{2d}$. 
 Let us define the following assumptions on $V$ and $c$:
\begin{enumerate}
\item[\textbf{(Av1)}] $V:\mathbb R^{d}\to \mathbb R$ is $\mathcal C^1$ and $V \text{ is lower bounded on } \mathbb R^d$. 

\item[\textbf{(Ac1)}] $c:\,  \mathbb R^{2d}\to  \mathrm M_{d}(\mathbb R)$ (the space of square matrices of size $d$ with real coefficients)  is continuous. In addition, there exist $ \eta>0$ and $L>0$, such that
$$ \forall v \in    \mathbb R^d, \vert x\vert \ge L: \ \frac 12\big[c(x,v)+c^T(x,v)\big]\ge \eta  \,I_{\mathbb R^d}.$$
Finally, for all $N>0$, 
$$
   \sup_{\vert x\vert \le N, v\in \mathbb R^d}\Vert c(x,v)\Vert_{\rm{H.S}}<+\infty,
$$
where $\Vert c(x,v)\Vert_{\rm{H.S}}$ is the Hilbert-Schmidt norm of matrix and where $c^T$ is the the transpose matrix of $c$.  
\item[\textbf{(A$\Sigma$)}] $\Sigma: \mathbb R^{2d} \to \mathbb R$  a $\mathcal C^\infty$ function,  uniformly Lipschitz over  $ \mathbb R^{2d}$, and  such that for some $\Sigma_0>0$ and $\Sigma_\infty>0$, 
$$\forall \mathsf x\in \mathbb R^{2d}, \ \Sigma_0\le \Sigma(\mathsf x)\le \Sigma_\infty.$$
\end{enumerate}
We also consider the less stringent assumption  than  \textbf{(Ac1)}:
\begin{enumerate}
\item[\textbf{(Ac0)}]
 $c:\,  \mathbb R^{2d}\to \mathrm M_{d}(\mathbb R)$ is continuous and 
 $$   \exists A\ge 0,  \forall x,v \in    \mathbb R^d: \ \frac 12\big[c(x,v)+c^T(x,v)\big]\ge -A  \, I_{\mathbb R^d}.$$
\end{enumerate}
 This will allow us to consider in particular  Langevin processes   with unbounded $v$-dependent damping  coefficient:   
\begin{equation}\label{eq.VDP2}
\text{for some }\ell_0>0, \,\,  \forall x,v\in \mathbb R^d, \ c(x,v)=\vert v \vert^{\ell_0} 
\end{equation}
and fast growing potential, in the sense that  there exist $\mathsf n_0>2$,  $r_0>0$, and $r>0$, for all $\vert x\vert \ge r_0$, 
\begin{equation}\label{eq.VDP3}
 V \text{ satisfies \textbf{(Av1)}}, \,  r^{-1}\vert x\vert ^{\mathsf n_0}\le V(x)\le r  \vert x\vert ^{\mathsf n_0} \text{ and }  r^{-1}\vert x\vert ^{\mathsf n_0}\le x\cdot \nabla V(x) . 
\end{equation}
Notice that when~\eqref{eq.VDP2} holds,  Assumption \textbf{(Ac1)} is not satisfied but \textbf{(Ac0)} is satisfied.
When $V$, $c$, and $\Sigma$ satisfy respectively \textbf{(Av1)},~\textbf{(Ac0)}, and \textbf{(A$\Sigma$)}, there is a unique weak solution to~\eqref{eq.hypo}    by \cite[Lemma 1.1]{Wu2001}, which is thus a (strong) Markov process. In the following we always assume that  \textbf{(A$\Sigma$)} holds.

As explained in the introduction above (see also~\cite{Ne17}), the process $(x_t,t\ge 0)$ is stuck during long periods of time in neighborhoods of the local minima of $V$. This is  due to energetic barriers this process has to cross to visit other regions of $\mathbb R^d$. Such neighborhoods are called metastable regions and are thus of the form 
\begin{equation}\label{eq.D}
\mathcal D= \mathsf O\times \mathbb R^d,
\end{equation}
where $\mathsf O\subset \mathbb R^d$. We are interested  in the existence of quasi-stationary distributions for the  processes~\eqref{eq.hypo} in domain $\mathcal D$ of this form.   
We assume that $\mathsf O$ is a $\mathcal C^2$ subdomain (not necessarily bounded) of $\mathbb R^d$, such that $\mathbb R^d\setminus \overline{\mathsf O}$ is non-empty.


  \subsubsection{The case when \textbf{(Av1)} and  \textbf{(Ac1)} hold}


  Let us define the following last assumptions on~$V$ and~$c$.  
\begin{enumerate}
\item[\textbf{(Av2)}]   There exists a $\mathcal C^1$ function $G:\mathbb R^d\to \mathbb R^d$ such that $G$ and $\nabla G$ are bounded over $\mathbb R^d$, and such that 
$$\nabla V(x)\cdot G(x)\to +\infty \text{ as } \vert x\vert \to +\infty.$$

\item[\textbf{(Ac2)}]  There exists some $\mathcal C^2$ lower bounded function $U:\mathbb R^d\to \mathbb R$ such that 
$$\sup_{x,v\in \mathbb R^d} \vert c^T(x,v)G(x)-\nabla U(x)\vert <+\infty.$$
\end{enumerate}
Some examples of functions $V$ and  $c$ satisfying \textbf{(Av1)}, \textbf{(Av2)}, \textbf{(Ac1)}, and \textbf{(Ac2)} are given in \cite[Remark 3.2]{Wu2001}, see also \cite[Remark 6.3]{guillinqsd}. 
  The Hamiltonian of the process~\eqref{eq.hypo} is, for  $x,v\in \mathbb R^d$,  
 $$\mathsf H(x,v)=V(x)+\frac12|v|^2.$$
Assume that   \textbf{(A$\Sigma$)}, \textbf{(Av1)}, \textbf{(Av2)}, \textbf{(Ac1)}, and \textbf{(Ac2)}  hold. 
Let us  introduce for $(x,v)\in \mathbb R^{2d}$, the modified Hamiltonian~\cite[Eq. (3.3)]{Wu2001}, 
 \begin{equation}\label{eq.F}
\mathsf F (x,v)=a \, \mathsf H  (x,v)+v\cdot (b \,   G(x)+\nabla \mathsf w(x)) +b \, U(x),
 \end{equation} 
where $G$, $U$ are as  \textbf{(Av2)} and \textbf{(Ac2)}, $a>0$, $b>0$, and $\mathsf  w:\mathbb R^d\to \mathbb R$ is a compactly supported $\mathcal C^2$ function.
Define, for all $x,v\in \mathbb R^d$:
 \begin{equation}\label{eq.Wkl}
\mathsf W  (x,v)=\exp\big[{\mathsf F(x,v) -\inf_{\mathbb R^{2d}} \mathsf F   }\big]\ge 1.
\end{equation} 
 The following  asymptotic upper bound on $\mathsf W_1$  gives a way to check   if a probability measure on $\nu$ on $\mathcal D\subset \mathbb R^{2d}$ satisfies $\nu(\mathsf W^{1/p})<+\infty$ (see the proof of  \cite[Lemma 2.6]{guillinqsd}). 
  \begin{lemma} 
Assume that $V$ satisfies \textbf{(Av1)} and \textbf{(Av2)}. Then $\lim_{x\to +\infty} V(x)=+\infty$. Let  $c$  be such that  \textbf{(Ac1)} and \textbf{(Ac2)} hold with $\lim_{\vert x\vert \to +\infty} U(x)/V(x)=0$. Then, for any $\epsilon>0$, there exists $R>0$ such that if  $\vert x\vert +\vert v\vert \ge R$,  $\mathsf W(x,v)\le e^{ a(1+\epsilon)\mathsf H(x,v)}$.  
\end{lemma}
We have the following result which ensures existence and uniqueness of 
  the quasi-stationary distribution of the process \eqref{eq.hypo} on $\mathcal D$  (see \eqref{eq.D})  when $V$ and $c$ satisfy \textbf{(Av1)}, \textbf{(Av2)}, \textbf{(Ac1)}, and \textbf{(Ac2)}.
  
  \begin{theorem}\label{thm:2}
Assume that $\Sigma$ satisfies~\textbf{(A$\Sigma$)} and that the functions  $V$ and $c$ satisfy \textbf{(Av1)}, \textbf{(Av2)}, \textbf{(Ac1)}, and \textbf{(Ac2)}. Then,  there exist parameters   $\mathsf w\in \mathcal C^2_c(\mathbb R^d,\mathbb R)$, $a>0$, and $b>0$  (see~\cite[Eq. (3.4) $\to$ Eq. (3.9)]{Wu2001} for explicit conditions on $\mathsf w$, $a$, and $b$) such that  Theorem~\ref{thm:main} is valid for the process~\eqref{eq.hypo} with  $\mathcal D=\mathsf O\times \mathbb R^d$  and with  the Lyapunov function $\mathsf W$ defined in~\eqref{eq.Wkl}. 
  \end{theorem}


In other words, for all $p>1$, under~\textbf{(A$\Sigma$)}, \textbf{(Av1)}, \textbf{(Av2)}, \textbf{(Ac1)}, and \textbf{(Ac2)} and   when $\mathcal D=\mathsf O\times \mathbb R^d$ (where $\mathsf O$ is as in Theorem \ref{thm:2}),  there exists a unique quasi-stationary distribution  in $\mathcal D$  for the process~\eqref{eq.hypo}   in    the space  $\mathscr M_p= \{\nu\in \mathcal M_1(\mathcal D), \nu(\mathsf W^{1/p})<+\infty\}$. In addition, Equation \eqref{eq.thm-maina} holds for all $\nu\in \mathscr M_p$.

  To prove Theorem \ref{thm:2} we use Theorem \ref{thm:main}. More precisely, we show that \textbf{(C1)}, \textbf{(C2)}, \textbf{(C4)},  and \textbf{(C5)} are satisfied and that, choosing well  the function $\mathsf w$ and  taking  $a>0$ and $b>0$  small enough (see~\cite[Eq. (3.4) $\to$ Eq. (3.9)]{Wu2001}), the function $\mathsf W$ defined in~\eqref{eq.Wkl} satisfies \textbf{(C3)}.

    When  $\mathsf O$ is bounded, $V$ is $\mathcal C^\infty$ on $\overline{\mathsf O}$, and both $\Sigma$   and   $c$ are   constant, the very recent work~\cite{LelRamRey2} (see also~\cite{LelRamRey,Ram}) establishes  the uniqueness of the quasi-stationary distribution in the whole space of measure over $\mathcal D$. These results were presented in Section~\ref{sec:Langevin 1} where the authors also considered non gradient vector fields (see indeed~\eqref{eq:Langevin 1}). Their approach is different from the one adopted in~\cite{guillinqsd,guillinqsd2}: it is based on a Feynman-Kac type formula for the killed semigroup, Harnack inequalities, and Gaussian upper bounds.  We finally also refer to~\cite{BenChaOcaVil}
   for the recent study of Hypoelliptic diffusions killed at the boundary of a bounded subdomain  of $\mathbb R^m$ with non characteristic boundary.


   \subsubsection{The case when \eqref{eq.VDP2} and \eqref{eq.VDP3} are satisfied}
  When \eqref{eq.VDP2} and \eqref{eq.VDP3} hold, and $  \ell_0<\mathsf n_0-2$, we are able to construct a bounded Lyapunov function for the process~\eqref{eq.hypo} (i.e. a bounded Lyapunov function $\mathsf W$ over $\mathbb R^{2d}$ satisfying \textbf{(C3)}). This implies the following theorem. 
  
  \begin{theorem}\label{thm:3}
Assume that $\Sigma$ satisfies~\textbf{(A$\Sigma$)} and that the functions  $V$ and $c$ satisfy    \eqref{eq.VDP2} and \eqref{eq.VDP3}. Assume in addition that 
$$
  \ell_0<\mathsf n_0-2.
$$
Then, Theorem \ref{thm:main} is valid with a bounded Lyapunov function. In particular:  there exists a unique quasi-stationary distribution in $\mathcal D=\mathsf O\times \mathbb R^d$  for the process~\eqref{eq.hypo}   in  the whole space of  probability measures $\mathcal M_1(\mathcal D)$ on $\mathcal D$, and in addition,  \eqref{eq.thm-maina} holds for all  $\nu\in \mathcal M_1(\mathcal D)$. 
  \end{theorem}


  \subsection{Langevin processes with singular potentials}
  \label{sec.Sec-appl2}
  
  In this section we consider the existence and uniqueness of a quasi-stationary distribution for the process \eqref{eq.hypo} when the potential $V$ is singular. For ease of exposition, we will just focus here on Lennard-Jones and Coulomb interactions. However much more general singular potentials are considered in~\cite{guillinqsd2}. 
  
  For  $d\ge 1$,  consider a system of $N\ge 2$ particles in $\mathbb R^d$ which cannot collide and let  
  $$x_t=(x^1_t,\ldots, x^N_t) \in  (\mathbb R^{d})^N \text{ and }  v_t=(v^1_t,\ldots, v^N_t) \in (\mathbb R^{d})^N,$$
denote respectively  the positions of the $N$ particles and   their velocities, at time $t\ge 0$. 
The natural space to consider the time evolution of the positions $(x_t,t\ge 0)$ and of the velocities  $(v_t,t\ge 0)$  of the $N$ particles is thus
\begin{equation}\label{eq.SS}
\mathcal S=\mathcal O \times   (\mathbb R^{d})^N,
\end{equation}
where, if $d=1$, 
$$\mathcal O= 
 \big \{x=(x^1,\ldots, x^N) \in    (\mathbb R)^N, \, x^1<x^2<\ldots<x^N \big \},$$
 and if $d\ge2$,
$$
   \mathcal O=  \big\{x=(x^1,\ldots, x^N)\in  (\mathbb R^{d})^N, \, x^i\neq x^j \text{ for all } i\neq j   \big \}. 
$$
  Notice that in both cases, $ \mathcal O$ is open, path connected, and unbounded.  In molecular dynamics,   the interatomic potential    of the system of $N$ particles is typically of the form, for $x=(x^1,\ldots,x^N)\in \overline{\mathcal O}$,
\begin{equation}\label{eq.VV}
 V(x)=\sum_{i=1}^N  V_{\textbf{{\rm c}}}(x^i)+\sum_{1\le i<j\le N}  V_{\textbf{{\rm I}}}(x^i-x^j)\in  \overline{\mathbb R},
\end{equation}
where $ V_{\textbf{{\rm c}}}:\mathbb R^d\to \mathbb R$    is the confining potential of the system  and $   V_{\textbf{{\rm I}}}:  {\mathsf  \Omega }\to  {\mathbb R}$  (where, if $d=1$, $ \mathsf  \Omega=\{y< 0\}$, and if $d\ge 2$, $ \mathsf  \Omega=\{y\neq 0\}$)   is  the potential energy modeling the interaction between two particles, the latter becoming infinite when  (and only when) $y\in \partial \mathsf \Omega=\{0\}$ (which prevents from collisions). 

Let   $(\Omega, \mathcal F, (\mathcal F_t)_{t\ge 0}, \mathbb P)$ be a  filtered probability space.  We assume that the evolution of the positions $(x_t,t\ge 0)$  and the velocities $(v_t,t\ge 0)$ of the  $N$ particles  on $\mathcal S$  is  described by the following  hypoelliptic stochastic differential equation     
  \begin{equation} \label{eq.hypo2}
\left\{
\begin{aligned}
 d  x_t &= v_t dt,\\
dv_t&=-\nabla   V(x_t)dt - c(x_t,v_t)   v_t dt +  \Sigma(x_t,v_t) \, dB_t,
\end{aligned}
\right.
\end{equation}
  where $c: (\mathbb R^{d})^N\times (\mathbb R^{d})^N \to \mathrm M_{Nd}(\mathbb R)$  is the friction matrix, $ ( B _t, t\ge 0  )$ is a standard $dN$-dimensional Brownian motion, and $\Sigma: (\mathbb R^{d})^N\times (\mathbb R^{d})^N \to \mathbb R$. We set $X_t=(x_t,v_t)$ for $t\ge 0$.  
  Throughout this section, we assume that  $\Sigma$ satisfies \textbf{(A$\Sigma$)} (with $Nd$ instead of $d$ there) and that $c$ is such that:
  \begin{enumerate}
    \item[\textbf{(Ac)}] $c:  (\mathbb R^{d})^N\times (\mathbb R^{d})^N \to \mathrm M_{Nd}(\mathbb R)$ is  a locally Lipschitz  function such that:   \begin{enumerate}
    \item[(i)] there exists  $c^*>0$,  $ \forall x,v \in    (\mathbb R^d)^N: \ \frac 12\big[c(x,v)+c^T(x,v)\big]\ge c^*  \,I_{(\mathbb R^d)^N}$, 
    \item[(ii)] $c$ is bounded, i.e. 
   $\sup_{x,v\in (\mathbb R^d)^N, k,\ell=1,\ldots,Nd}| c_{k,\ell}(x,v)|<+\infty$. 
      \end{enumerate}
\end{enumerate}
In the two next sections, we consider  existence and uniqueness of a quasi-stationary distribution for  the process~\eqref{eq.hypo2} when  $ V_{\textbf{{\rm I}}}$ 
   is the Lennard-Jones   potential or the Coulomb potential.


\subsubsection{Quasi-stationary distributions for the Lennard-Jones potential (for any $d\ge 1$) and the Coulomb potential (when $d\ge 3$)} 
We recall that the Lennard-Jones potential is defined by, for   $y\in \overline{\mathsf \Omega}$, 
  $$ V_{LJ}(y) \!=\!\left\{
 \begin{aligned}
 \frac{b }{\vert y\vert^{{12}}}-\frac{c}{\vert y\vert^6}  &\text{ if } y\in \mathsf  \Omega\\
 +\infty  &\text{ if } y= 0,
 \end{aligned}
 \right. $$ 
 where $b,c>0$. 
When $d\ge 3$,   the Coulomb potential is the function defined on  $\overline{\mathsf \Omega}$ by:
   $$
    V_{co}(y) \!=\!\left\{
 \begin{aligned}
 \frac{e}{\vert y\vert^{d-2}}  &\text{ if } y\in \mathsf  \Omega\\
 +\infty  &\text{ if } y= 0,
 \end{aligned}
 \right. $$ 
 where $e>0$.         
  In this section, the interaction potential $V_{\textbf{{\rm I}}}$ is either $  V_{co}$ or $V_{LJ}$.   The Coulomb potential when $d=2$ is treated in the next section. We assume that  $ V_{\textbf{{\rm c}}}$ is $\mathcal C^2$ over $\mathbb R^d$ and to simplify we assume that for some $r>0$, 
\begin{equation}\label{eq.Vcc}
\forall y\in \mathbb R^d, \vert y\vert \ge r,  \   V_{\textbf{{\rm c}}}(y)= A \vert y\vert^{\alpha_1},
\end{equation}
  where $ A>0$ and $\alpha_1>1$. 
We refer to \cite{guillinqsd2} for more general confining potentials. 
With such confining potentials and interaction potentials, by \cite[Proposition 2.3]{guillinqsd2}, for any $\mathsf x\in \mathcal S$, there exists a   unique pathwise solution $( X _t=(x_t,v_t), t\ge 0  )$  of \eqref{eq.hypo2} with $X_0=\mathsf  x$, which is moreover non-explosive and remains in $\mathcal S$ for all $t\ge 0$. 
 
 Recall that the  Hamiltonian of the process~\eqref{eq.hypo}   is, for  $(x,v)\in  \mathcal S$,  
  \begin{equation}\label{eq.Ha}
 \mathsf H(x,v)=  V(x)+\frac12|v|^2,
  \end{equation}
  where $V$ is given by \eqref{eq.VV}. 
  We have that $ \mathsf H(x,v)\to +\infty$ (for $(x,v)\in \mathcal S$, see \eqref{eq.SS}) if and only if $x\to \partial \mathcal O\cup\{\infty\}$ or $v\to +\infty$.

  Let us now consider 
$$\eta_1\in (0,1].$$
If $\alpha_1\in (1,2)$, we assume in addition  that $\eta_1>(2-\alpha_1)/\alpha_1$ ($\in (0,1)$). 
We have the following result (see \cite[Theorem 2.4]{guillinqsd2}).


\begin{theorem}\label{thm:4}
Assume that   \textbf{(Ac)} and \textbf{(A$\Sigma$)} are satisfied. Assume also  that $V_{\textbf{{\rm I}}}\in \{  V_{co},V_{LJ}\}$  and that $  V_{\textbf{{\rm c}}}$ is a $\mathcal C^2$ potential such that \eqref{eq.Vcc} holds. 
 Let  $\mathsf O$ be a   subdomain  of $\mathcal O$ such that $\mathcal O\setminus \overline{\mathsf O}$ is nonempty  and $\partial{\mathsf O}\cap  \mathcal O$ is $\mathcal C^2$. Set   $\mathcal D=\mathsf O\times( \mathbb R^{d})^N\subset \mathcal S$ (see \eqref{eq.SS}). Choose $\eta_1$ as above. 
Then, Theorem \ref{thm:main} is valid for the process \eqref{eq.hypo2} on $\mathcal D$   
with a    continuous and unbounded Lyapunov functional $\mathsf W :\mathcal S\to [1,+\infty)$ which satisfies, for some $ m>0$,  $\mathsf W \le  \exp\big[m  \mathsf H    ^{\eta_1}\big]$ on $\mathcal S$. 
\end{theorem}
In other words,  Theorem \ref{thm:4} states that on $\mathcal D$,  there exists a unique quasi-stationary distribution for the  process \eqref{eq.hypo2}  in the space $\mathscr M_p$ ($\forall p>1$). In addition, Equation \eqref{eq.thm-maina} holds for all $\nu\in \mathscr M_p$. 

We refer to   \cite[Proposition 2.10]{guillinqsd2} for the explicit construction of $\mathsf W $, which is inspired by the previous works \cite{Wu2001}  and  \cite{herzog}. 
Notice that $\mathsf O$ is not necessarily bounded in   Theorem~\ref{thm:4}, and its closure may contain  singularities of $ V$, namely some subset of $\partial \mathcal O$. 
  

\subsubsection{Quasi-stationary distributions for  the Coulomb potential when $d=1,2$}
 In this section, we consider the interaction potential    $ V_{\textbf{{\rm I}}}$ defined by 
  \begin{equation}\label{de.Cou}
\text{ for all $y\in \mathbb R^d$: } V_{\textbf{{\rm I}}}(y)= -\log \vert y\vert \text{ if  } y\neq 0, \text{ else }   V_{\textbf{{\rm I}}}(y)=+\infty.
\end{equation}
When $d=2$, $ V_{\textbf{{\rm I}}}$ is the   Coulomb potential, and when $d=1$, $ V_{\textbf{{\rm I}}}$ corresponds to a  log singularity pairwise potential. We also assume that the confining potential $ V_{\textbf{{\rm c}}}:  \mathbb R^d  \to \mathbb R$   satisfies: $ V_{\textbf{{\rm c}}}\in \mathcal C^2 (\mathbb R^d,\mathbb R)$ and for some  $r>0$, it holds:
\begin{equation}\label{eq.Vcc2}
\forall y\in \mathbb R^d, \vert y\vert \ge r,  \  V_{\textbf{{\rm c}}}(y)=  A \vert y\vert^{\alpha_2},
\end{equation}
  where $A>0$ and $\alpha_2\ge0$ (we refer to \cite{guillinqsd2} for more general confining potentials).

  With such  an interaction potential and confining potentials, by \cite[Proposition 3.1]{guillinqsd2}, for any $\mathsf x\in \mathcal S$ (see \eqref{eq.SS}), there exists a   unique pathwise solution $( X _t=(x_t,v_t), t\ge 0  )$  of \eqref{eq.hypo2} with $X_0=\mathsf  x$, which is moreover non-explosive and remains in $\mathcal S$ (see \eqref{eq.SS}) for all $t\ge 0$. 
  
  \begin{theorem}\label{thm:5} Let $d\in \{1,2\}$. 
Assume that   \textbf{(Ac)} and \textbf{(A$\Sigma$)} are satisfied.  Assume also  that $V_{\textbf{{\rm I}}}$ is given by \eqref{de.Cou}   and that $  V_{\textbf{{\rm c}}}$ is a $\mathcal C^2$ potential such that \eqref{eq.Vcc2} holds.  Let  $\mathsf O$ be a   subdomain  of $\mathcal O$ such that $\mathcal O\setminus \overline{\mathsf O}$ is nonempty  and $\partial{\mathsf O}\cap  \mathcal O$ is $\mathcal C^2$. Set   $\mathcal D=\mathsf O\times( \mathbb R^{d})^N\subset \mathcal S$ (see \eqref{eq.SS}). Take $\eta_2\in (0,1]$. 
 Then, Theorem \ref{thm:main} is valid for the process \eqref{eq.hypo2} on $\mathcal D$   
 with a    continuous and unbounded Lyapunov functional $\mathsf W :\mathcal S\to [1,+\infty)$ which satisfies, for some $ m>0$,  $\mathsf W \le  \exp\big[m  \mathsf H    ^{\eta_2}\big]$ on $\mathcal S$. 
\end{theorem}

We refer to   \cite[Proposition 3.3]{guillinqsd2} for the explicit construction of $\mathsf W $, which is inspired by~\cite{lu2019geometric}. 
Let us emphasize  that   there is no restriction in Theorem \ref{thm:5}  on $\eta_2$ (i.e. one can choose any $\eta_2$ in $(0,1]$). We mention that in \cite{guillinqsd2}, we also deduce that large deviation principles hold for the  non killed process  (i.e. the process $(X_t,t\ge 0)$ on $\mathcal S$) with Lennard-Jones and Coulomb potential interactions, as well as the exponential convergence of its law towards its invariant measure. The existence and uniqueness of a quasi-stationary distribution for elliptic diffusions with Lennard-Jones type interactions is also studied there.

\section{Numerical approximation with particles}
\label{sec:numeric}

This section is based on \cite{JournelMonmarche}. The goal is to introduce a numerical scheme  whose aim is to sample the  QSD of a killed diffusion, and to prove quantitative long-time convergence rates for this algorithm.

\subsection{The problem}

Given a Markov process $X$ killed at time $\tau$, we introduce a numerical scheme in order to sample its quasi-stationary distribution $\nu_{\ast}$, or the law of $X_t$ conditioned on survival. With great generality, one may introduce a particle approximation of the process conditioned on survival, called a Fleming-Viot process in most of the literature. It consists in $N\in\mathbb N$ independent processes $(X^i)$, having the same dynamic as $X$, but when one of them dies, instead of staying at the cemetery point, it chooses uniformly at random another one of the processes, and branches onto it. This gives rise to a mean-field interacting particle system, with Markov semi-group $(P_{N,t})_{t\geqslant 0}$. In the particular case of a softly-killed elliptic SDE on the $d$-dimensional torus:
\begin{equation}\label{diffontorus}
\dd X_t = b(X_t)\dd t + \dd B_t,
\end{equation}
where $b\in\mathcal C^1(\T^d)$ and $B$ is a $d$-dimensional Brownian motion, we propose a numerical scheme of this Fleming-Viot process.

Given a random variable $E\in\R_+$ with law $e^{-t}\dd t$, we define the soft killing as:
\begin{equation}
 \tau = \inf\left\{t\geqslant0,E\leqslant\int_0^t\lambda(X_s)\dd s\right\},
\end{equation}
where $\lambda:\T^d\to\R_+$ is a $L_\lambda$-Lipschitz death rate, for some $L_\lambda>0$. As before, write:
\[\mu_t=\mathcal Law(X_t|\tau>t),\]the law of the process conditioned on survival.

First, we introduce a discretized version of the killed diffusion. We define a Markov kernel $K^{\gamma}$ on $\T^d$ for $\gamma$ small enough as follow: fix $x\in\T^d$, then \[K^{\gamma}(x,\cdot)=\mathcal N(x+\gamma b(x),\gamma I).\] We define a Markov chain $( X_n)$ by $X_{n+1}\sim K^{\gamma}(X_n,\dd x)$. Given a sequence of independent and uniform on $[0,1]$ random variable $(U_n)$, we define the death time of the chain $(X_n)$ as:
\begin{equation}
\bar \tau = \inf\left\{n\geqslant 0, U_n\leqslant p(X_n)\right\},
\end{equation}
where $p(x)=1-e^{-\gamma \lambda(x)}$. Write \[\eta_n=\mathcal{L}aw(X_n|\tau>n).\] The discretized killing as been defined as such in order to have for all $t\geqslant 0$:
\[\eta_{\lfloor t/\gamma\rfloor} \rightarrow \mu_t,\]as $\gamma\rightarrow 0$. This Markov chain also admits a QSD $\nu_\gamma$.

Now, we want to define two processes: a Markov process $(Y_n)$ whose law is $\eta_n$, defined for all time, and a particle approximation of this Markov chain $(Y_n)$. To this end, we define first, given some probability measure $\mu\in\mathcal M^1(\T^d)$, a Markov kernel $Q_\mu$ such that $Q_\mu(x,\cdot)$ corresponds to the law of some random variable having done a step of discretized diffusion, and upon death, it does a step of discretized diffusion, conditioned on survival, with initial condition $\mu$. More precisely, we construct a random variable of law $Q_\mu(x,\cdot)$ as follow:
\begin{enumerate}
\item Draw  $X_{0} \sim K^{\gamma}(x,\cdot)$ and $U_0\sim \mathcal U([0,1])$.
\item If $U_0 \geqslant p(X_{0})$, set $X = X_{0}$ in $\mathbb T^d$ (in that case, we say the particle has moved from $x$ to $X_0$ without dying).
\item If $U_0 <p(X_{0})$ then set $i=1$ and, while $X$ is not defined, do:
\begin{enumerate}
\item Draw a new $X_{i}'$ distributed according to $\mu$, a new $X_{i} \sim \mathcal N\po X_{i}' + \gamma b(X_{i}') , \gamma I_d \pf$ and a new $U_i\sim \mathcal U([0,1])$.
\item  If $U_i \geqslant p(X_{i})$, set $X = X_{i}$ in $\mathbb T^d$  (in that case, we say the particle has died, resurrected at $X_{i}'$, moved to $X_{i}$ and survived).
\item  If $U_i < p(X_{i})$, set $i\leftarrow i+1$ (in that case, we say the particle has died, resurrected at $X_{i}'$, moved to $X_{i}$ and died again) and go back to step (a).
\end{enumerate}
\end{enumerate}

We introduce the following non-linear process, solution to the problem:
\begin{equation}\label{nonlineaire}
\left\{\begin{matrix} Y_{n+1}\sim Q_{\tilde \eta_n}(Y_n,\cdot) \\ \tilde \eta_n=\mathcal Law(Y_n) \end{matrix} \right. .
\end{equation}
This non-linear Markov chain does at each time $n$ a step of discretized diffusion conditioned on survival. As said, the reason to introduce this Markov chain is the following property:
\begin{proposition}\label{Prop-EgaliteLoi}
For all $n\in\mathbb N$
\[\tilde \eta_n \ = \ \mathcal Law\po X_n \ | \ \tau \ > \ n\pf\,.\]
\end{proposition}

Now we can define the numerical scheme of the Fleming-Viot process. For $\textbf x=(x_1,\dots,x_N)$, we write:
\begin{equation}
\pi(\textbf x) \ := \ \frac{1}{N}\sum_{i=1}^N \delta_{x_i}\ \in \mathcal P( \mathbb T^d)
\end{equation}
the empirical measure and:
\begin{equation}
R_{\gamma,N}(\textbf x,\cdot) = Q_{\pi(x)}(x_1,\cdot)\otimes \dots \otimes Q_{\pi(x)}(x_N,\cdot).
\end{equation}
We then define recursively the process by $\textbf{X}_{n+1}\sim R_{\gamma,N}(\textbf X_n,\cdot)$. This can be understood as follow: each particles is an Euler-Maruyama scheme of the diffusion, and if one of them dies at time $n$, then it does a step of $K^\gamma$, starting from $\pi(\textbf X_{n-1})$, conditioned on survival. As $N$ goes to infinity, if the initial condition condition of this particle system is $N$ iid random variable of law $\eta_0$, then $\pi(\textbf X_n)$ is expected to be close to the common law of the $X^i$'s. This can then be shown to be close to $\eta_n$. This is the propagation of chaos phenomenon.

\subsection{Main results}

For this Markov chain, we proved two theorem. The first one is about the long-time behavior of the process:

\begin{theorem}\label{Theorem-Inter}
There exists $c_0,\gamma_0>0$ that depends only on the drift $b$ and the dimension $d$ such that, if $\lambda$ is Lipschitz with a constant $L_\lambda$ such that
\begin{eqnarray}\label{Eq-ConditionPerturb}
L_\lambda e^{\gamma \|\lambda\|_\infty} & < & c_0\,,
\end{eqnarray}
then there exist a distance $\rho$  equivalent to $|\cdot|$ on $\T^d$, such that $\rho_N(x,y)=\sum \rho(x_i,y_i)$ is a distance on $\T^{dN}$, and $\kappa>0$, such that for all $\gamma\in (0,\gamma_0]$ 	$N\in\N$, and all $\textbf x,\textbf y\in\T^{dN}$, there exists a coupling $(\textbf X,\textbf Y)$ of $\delta_{\textbf x}R_{N,\gamma}$ and $\delta_{\textbf y}R_{N,\gamma}$ such that:
\begin{align*}
\mathbb E\left( \rho_N\left(\textbf X,\textbf Y\right) \right)\leqslant (1-\gamma\kappa)\rho_N(\textbf x,\textbf y)\,.
\end{align*}
This yields the existence of a stationary measure $\mu_{\infty,N,\gamma}$, and exponential convergence of the Markov chain towards this probability measure.
\end{theorem}

To construct the coupling of theorem~\ref{Theorem-Inter}, we first construct a coupling of $Q_\mu(x,\cdot)$ and $Q_\nu(x,\cdot)$, for any $x,y\in \T$, $\mu,\nu\in\mathcal P(\T^{d})$. This coupling is done as follow: thanks to \cite{EberleMajka}, we are able couple the step of discretized diffusion to get a coupling $(X',Y')$ of $\delta_xK^\gamma$ and $\delta_yK^\gamma$. The probability that one of them dies and not the other is less then $L_\lambda \mathbb E(|X'-Y'|)$. If they die together, then we construct again thanks to \cite{EberleMajka} a sequence of coupling of $\mu K^\gamma$ and $\nu K^\gamma$.

Replacing in this coupling $\mu$ and $\nu$ by $\pi(\textbf x)$ and $\pi(\textbf y)$ yields a coupling of $R_{\gamma,N}(x,\cdot)$ and $R_{\gamma,N}(y,\cdot)$.

The second theorem is about the convergence of the empirical measure of the process towards the quasi-stationary distribution of the process~\eqref{diffontorus} killed at time $\tau$, in the limit $\gamma\rightarrow 0$ and $n,N\rightarrow\infty$. 

\begin{theorem}\label{Theorem-Main}
Under the hypothesis of theorem \ref{Theorem-Inter}, there exists $C>0$ such that for all $N\in\N$, $\gamma\in (0,\gamma_0]$, $t\geqslant 0$ and $\mu_0\in\mathcal P(\mathbb T^{dN})$, if $(\textbf X_k)_{k\in\N}$ is a Markov chain with initial distribution $\mu_0$ and transition kernel $R_{N,\gamma}$, 
\[\mathbb E \co \mathcal W_1\po \pi( \textbf X_{\lfloor t/\gamma\rfloor}), \nu \pf\cf \ \leqslant \ C \po  \sqrt \gamma + \alpha(N) + e^{-\kappa t}\pf\,,\]
where
\[\alpha(N) \ = \ \left\{\begin{array}{ll}
N^{-1/2} & \text{if }d=1\,,\\
N^{-1/2}\ln(1+N) & \text{if }d=2\,,\\
N^{-1/d} & \text{if }d>2\,.
\end{array}\right.\]
\end{theorem}  

The second theorem is a consequence of the first one. The rate of convergence in $\gamma$, $N$ and $t$ would be the same in the case of $N$ independent diffusion, hence are optimal. The rate of convergence in $N$ comes from \cite{FournierGuillin}.

\begin{figure}[h]
\begin{center}
$
\xymatrix{
    R_{N,\gamma}^m \ar[rrr] \ar[ddd] \ar[dr] &&& P_{N,t} \ar[dr] \ar[ddd] \hole \\
    & \eta_m \ar[rrr] \ar[ddd] &&& \mu_t \ar[ddd] \\
		&&&&\\
    \mu_{\infty,N,\gamma} \ar[rrr] \hole \ar[dr] &&& \bar{\mu}_{\infty,N} \ar[rd] \\
    & \nu_{\gamma} \ar[rrr] &&& \nu \\
  }
\xymatrix{
 \ar[dd]_{t\rightarrow\infty} \ar[dr]^{N\rightarrow\infty} \ar[rr]^{\gamma \rightarrow 0} && \\ && \\ &&
 }
$
\end{center}
\caption{Summary of the results}
\label{cubeconv}
\end{figure}
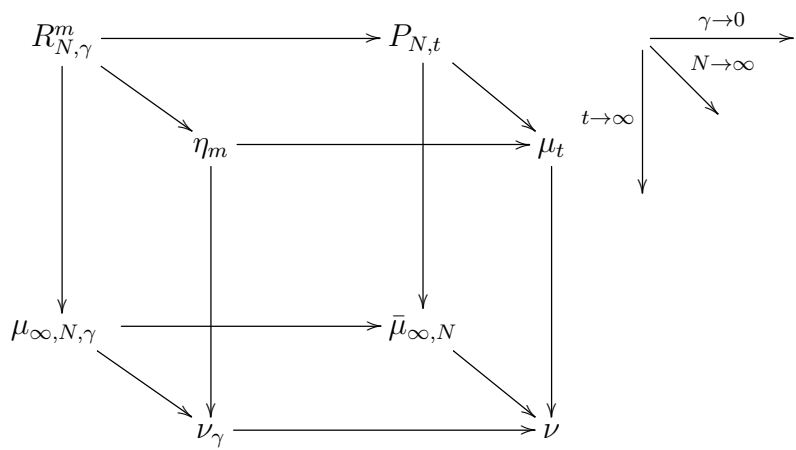

To prove propagation of chaos, we construct a coupling of the discretized Fleming-Viot process, with a system of $N$ independent non-linear process defined in equation~\eqref{nonlineaire}. This coupling is done is a similare fashion as for the previous coupling of $R_{\gamma,N}(x,\cdot)$ and $R_{\gamma,N}(x,\cdot)$. For the limit $\gamma\rightarrow 0$, the coupling of the process is made with a continuous time Fleming-Viot process.  Actually, thanks to those couplings, we are able to prove similar convergence when only one or two parameters among $t$, $\gamma$, or $N$ are sent to their limit. We can summarize our results thanks to the diagram of figure~\ref{cubeconv}.

\bibliography{biblio-acte}
\bibliographystyle{plain}
\end{document}